\documentclass[10pt,psamsfonts,french]{amsart}

\usepackage[french]{babel}
\usepackage[T1]{fontenc}    

\usepackage{amsmath}
\usepackage{amssymb}
\usepackage{amscd}
\usepackage{amsfonts}
\usepackage{amsbsy}
\usepackage{graphicx}
\usepackage[dvips]{psfrag}
\usepackage{url}
\usepackage{epsfig}
\usepackage{latexsym}
\usepackage[mathscr]{eucal}
\usepackage{graphics}

\begin{document}

\title[THE PARADOX OF VITO VOLTERRA'S PREDATOR-PREY MODEL]
{THE PARADOX OF VITO VOLTERRA'S \\ 
PREDATOR-PREY MODEL}

\author[J.M. Ginoux]
{Jean-Marc Ginoux$^{1,2}$}

\address{$^1$ Laboratoire LSIS, CNRS, UMR 7296, Universit\'{e} de Toulon,
BP 20132, F-83957 La Garde cedex, France}
\address{$^2$ Archives Henri Poincar\'{e}, Université de de Lorraine, CNRS, UMR 7117, BP 454, F-54001 Nancy, Cedex, France}

\email{ginoux@univ-tln.fr}

\maketitle

\begin{center}
\textit{This article is dedicated to the late Giorgio Israel.}
\end{center}

\begin{abstract}
The aim of this article is to propose on the one hand a brief
history of modeling starting from the works of Fibonacci, Robert Malthus,
Pierre Francis Verhulst and then Vito Volterra and, on the other hand, to
present the main hypotheses of the very famous but very little known predator-prey
model elaborated in the 1920s by Volterra in order to solve a problem
posed by his son-in-law, Umberto D'Ancona. It is thus shown that, contrary
to a widely-held notion, Volterra's model is realistic and his seminal work laid
the groundwork for modern population dynamics and mathematical ecology,
including seasonality, migration, pollution and more.
\end{abstract}

\section{A short history of modeling}

\subsection{The Malthusian model.}
\label{sec1}

\begin{figure}[htbp]
\centerline{\includegraphics[width = 8cm,height = 8cm]{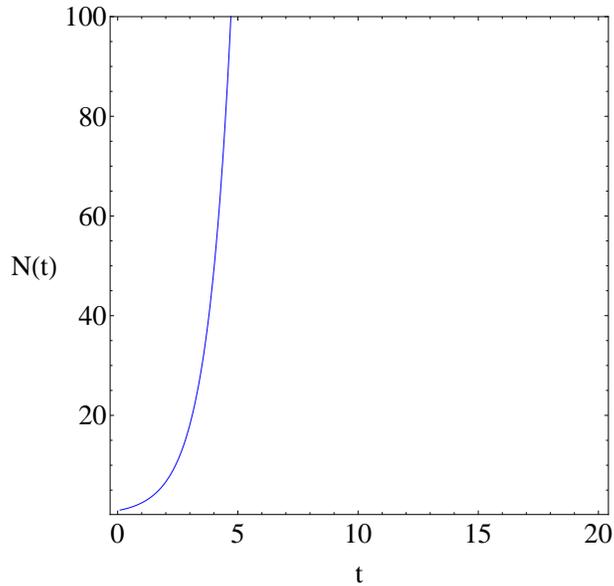}}
\caption{Malthus's model of exponential growth.}
\label{fig1}
\end{figure}

If the first scientific view of population growth
seems to be that of Leonardo Fibonacci \cite{fib}, also called Leonardo of Pisa, whose
famous sequence of numbers was presented in his Liber abaci (1202) as a solution
to a population growth problem, the modern foundations of \textit{population dynamics}\footnote{According to Frontier and Pichod-Viale \cite{fron} the correct terminology should be population ``kinetics'', since the interaction between species cannot be represented by forces.}
clearly date from Thomas Robert Malthus \cite{malt}. Considering an ``ideal'' population\footnote{A population is defined as the set of individuals of the same species living on the same territory
and able to reproduce among themselves.}
consisting of a single homogeneous \textit{animal species}, that is, neglecting the variations
in age, size and any periodicity for birth or mortality, and which lives alone in
an \textit{invariable environment} or coexists with other species without any direct or
indirect influence, he founded in 1798, with his celebrated claim ``Population, when
unchecked, increases in a geometrical ratio'', the paradigm of exponential growth.
This consists in assuming that the increase of the number $N\left( t \right)$ of individuals of this
population, during a short interval of time, is proportional to $N\left( t \right)$. This translates
to the following differential equation:

\begin{equation}
\label{eq1}
\frac{dN\left( t \right)}{dt} = \varepsilon N\left( t \right)
\end{equation}

where $\varepsilon$ is a constant factor of proportionality that represents the \textit{growth coefficient}
or \textit{growth rate}. By integrating (\ref{eq1}) we obtain the \textit{law of exponential growth} or \textit{law of
Malthusian growth} (see Fig. \ref{fig1}). This law, which does not take into account the limits
imposed by the \textit{environment} on growth and which is in disagreement with the actual
facts, had a profound influence on Charles Darwin's work on natural selection.
Indeed, Darwin \cite{darw} founded the idea of ``survival of the fittest'' on the impossibility
of an indefinite population growth. He illustrated this impossibility by a superb parabola describing the descendants of a pair of elephants that, under optimal
conditions, would cover the surface of the Earth in a few centuries\footnote{An example of this impossibility is also shown in the film by V.A. Kostitzin and J. Painlevé titled: ``Images mathématiques de la lutte pour la vie'', 1937, Médiatèque du Palais de la
Découverte, Paris.}. However, in laboratory experiments, the predictions of the \textit{Malthusian law} remain correct on
small numbers, while there is divergence for high values of the population.
Thus, we are led to conclude that the \textit{exponential law} remains valid as long as
the density of the population does not saturate the \textit{environment}.

\subsection{P.F. Verhulst's model.}

\begin{figure}[htbp]
\centerline{\includegraphics[width = 8cm,height = 8cm]{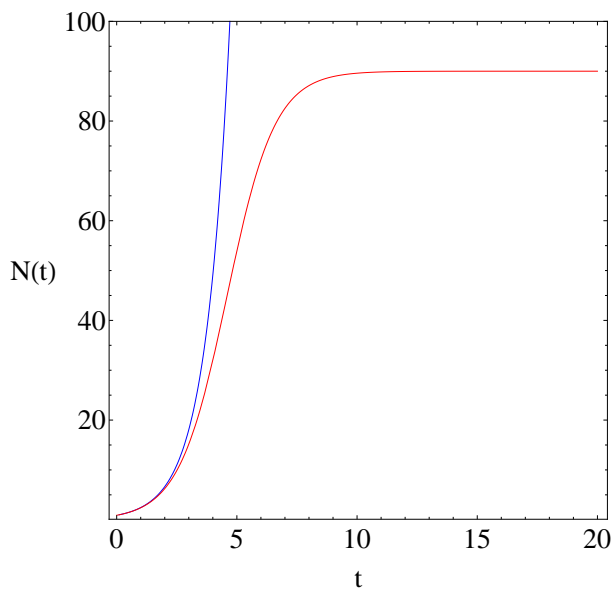}}
\caption{Verhulst's model of logistic growth (in red).}
\label{fig2}
\end{figure}

It was on the basis of these considerations that the
Belgian biologist Pierre-François Verhulst \cite{verh} proposed in 1837 a model that took
into account the limitation imposed by the increasing population size:

\begin{equation}
\label{eq2}
\frac{dN\left( t \right)}{dt} =  \varepsilon N\left( t \right) - \lambda N^2\left( t \right) = \varepsilon N\left( t \right)\left( 1 -
\frac{1}{K} N\left( t \right) \right)
\end{equation}

where $\varepsilon$ represents the \textit{growth rate}. The second coefficient $\lambda = \varepsilon / K $ originates in a
``mechanistic'' interpretation of the phenomenon. Indeed, it is assumed that growth
is limited by a kind of interior ``friction'' within the population, that is, the resources
remaining the same, the higher the number of individuals, the more difficult it is for
them to feed themselves and then to grow. This is a struggle between individuals
for the existence, and so an intra specific competition for food. The factor $K$, called
\textit{carrying capacity}, corresponds to the \textit{capacity} of the \textit{environment} to support the
population growth and represents the \textit{population limit} beyond which it can no longer
grow. This law, which Verhulst called \textit{logistic equation}, is radically different from
Malthus's, since it imposes a limiting value on the population (see Fig. \ref{fig2}, where exponential growth is represented in blue and logistic growth in red). It has been successfully applied to many real life situations, such as population growth in the
United States between 1790 and 1950 (see Pearl and Reed \cite{pearl}) or in experiments
conducted by the Russian biologist Georgii Frantsevich Gause \cite{gau} on the growth of
a protozoan, \textit{Paramecium caudatum}. Independently of these two archetypes, other
growth models have been developed: let us mention, for instance, the model by B.
Gompertz \cite{gom} aimed at evaluating the growth rate of a tumour.

\subsection{Volterra's ``predator-prey'' model.}

In the first half of the twentieth century,
the study of the dynamics of several interacting species developed considerably.
It was at this time, called ``the golden age of theoretical ecology'' \cite{scud}, that the
first models based on competition-type behaviours and predator-prey relationships
were developed. The paternity of the first model used to transcribe this kind of
interactions was the subject of a quarrel between Alfred J. Lotka \cite{lotka} and Vito
Volterra\footnote{For more about the dispute over priority, see \cite{guerra}, pp. 146-147.}, which is described in the work of Giorgio Israel \cite{isra91, isra93, isra96, isra02}. It has been established that the merit of the development of the first model of predator-prey
type applied to a biological situation is Volterra's. Indeed, Volterra's interest in problems of equilibria between animal species in ecosystems was prompted by
his son-in-law, zoologist Umberto D'Ancona, who for some years had been dealing
with statistics on fishing in the northern Adriatic Sea. These data involved the percentage
of predatory fish (Selachians) caught in three Italian ports, Trieste, Fiume
(now Rijeka) and Venice, during the period 1905-1923. They proved that during
the period 1915-1920, when fishing was less intense because of the war, there had
been a relative increase in the Selachians. According to D'Ancona's hypothesis, fishing disturbed the natural balance between species. It favoured a relative increase in ``prey'' species, that is, fish that feed only on plankton, and a decrease in
``predatory'' species, that is, fish that feed on other fish. The decline in fishing due
to World War I had thus restored, at least in part, the natural balance. D'Ancona
turned to Volterra, asking him to find a mathematical proof of his hypothesis. In
1926 Volterra published a reply in Italian \cite{vol26} which took the form of the famous
``predator-prey'' model, reproduced below. He then presented his results in English
and, in more detailed and complete form, in French \cite{vol27,vol28, vol31}. In the first
chapter of his ``Leçons sur la théorie mathématique de la lutte pour la vie'' \cite{vol31},
Volterra studies the coexistence of ``two species, one of which devours the other''.
Considering two species, the first, the prey $N_1\left( t \right)$, would undergo a Malthusian
growth if it were alone. The second, the predator $N_2\left( t \right)$, feeds exclusively on the
first and in the absence of prey will gradually wear out and disappear. The formulation
of the equation representing the predation is based on the \textit{méthode des
rencontres} (``method of encounters'') and on the \textit{hypothèse des équivalents} (``hypothesis
of equivalents'') due to Volterra \cite{vol31}. The former assumes that for predation
to occur between a predatory species and a prey species, it is first necessary to have
encounters between these two species and that the number of encounters between
them is proportional to the product of the number of individuals composing them,
that is, $N_1\left( t \right) N_2 \left( t \right)$, the coefficient of proportionality being equal to the probability
of an \textit{encounter}. The second hypothesis consists in assuming that ``there is a constant
ratio between the disappearances and appearances of individuals caused by
the encounters'', that is, that predation of the preys is \textit{equivalent} to increase of the
predators. At the beginning, Volterra considers this increase as \textit{immediate}\footnote{This means that predation is immediately transcribed in terms of growth of the predator species, whereas its effect naturally occurs with some delay. It will be seen below that Volterra also took this delay into account.}. This led him to the system:

\begin{equation}
\label{eq3}
\left\{
\begin{aligned}
\dfrac{dN_1}{dt} & = \varepsilon_1 N_1 - \gamma_1 N_1 N_2 = N_1 \left( \varepsilon_1 - \gamma_1 N_2  \right) \hfill \\
\dfrac{dN_2}{dt} & = - \varepsilon_2 N_2 + \gamma_2 N_1 N_2 = - N_2 \left( \varepsilon_2 - \gamma_2 N_1  \right)
\end{aligned}
\right.
\end{equation}

where $\varepsilon_1$ represents the prey's growth rate in the absence of the predator;
$\gamma_1$ the predation rate of the predator on the prey; $\varepsilon_1$ the predator's mortality rate in the
absence of prey; and $\gamma_2$ the growth rate of the predator due to its predation. From
this model, Volterra was able to state the \textit{law of the disturbance of the averages}:

\begin{quote}
If an attempt is made to destroy the individuals of the two species
uniformly and in proportion to their number, the average of the
number of individuals of the species that is eaten increases and that
of the individuals of the species feeding upon the other diminishes\\ \vphantom{} \hfill \cite[p. 20]{vol28}
\end{quote}

To establish this result, Volterra assumes that, for a time interval $dt$, we destroy $\alpha \lambda N_1 dt$ preys and $\beta \lambda N_2 dt$ predators. He then proves that the average values of prey and predators, which were previously equal to $\varepsilon_2 / \gamma_2$ and $\varepsilon_1 / \gamma_1$, respectively,
become, after this destruction, $( \varepsilon_2 + \beta \lambda )  / \gamma_2$ and $( \varepsilon_1 - \alpha \lambda)  / \gamma_1$. The phenomenon
observed by D'Ancona is thus explained: the increase in the number of predators and the reduction in the number of prey resulted from the decline in fishing, which before the war had changed the natural balance of this ``biological association''. A
decrease in fishing favours the more voracious species at the expense of the other.
This deterministic model, which constitutes the archetype of the trophic network,
aims at transcribing different types of animal behaviour into mathematical functions.
Two types of behaviour are represented: those related to \textit{increase} and those
related to \textit{decrease}. \textit{Natality} and \textit{predation} are related to \textit{increase}, whereas \textit{natural
mortality} and \textit{mortality by predation} correspond to a decrease in the number of
individuals. Each of these behaviours has a mathematical form, called a \textit{functional
response}.

\subsection{Different types of functional responses.}

Since the mid-1920s, these \textit{functional responses} have been the subject of numerous studies and developments aimed
at making the representation of animal behaviour by a mathematical function more
realistic. Natural growth, that is by natality, of the prey, represented by a functional
response of Malthus type \cite{malt} was later modified by Verhulst \cite{verh} to account for
its being bounded. The decrease by natural mortality was initially considered in a
way symmetrical to that of natural growth, that is, by substituting in the equations
(\ref{eq1}) and (\ref{eq2}) the growth rate $\varepsilon$ with a natural mortality rate $-\varepsilon$. It is important to
emphasise that this functional response, sometimes called \textit{closure relation}, was later
the subject of special studies aimed at transcribing the specific behaviour of certain
species, for instance cannibalism. All these functional responses are summarized in
Table 1.

\begin{table}[htbp]
\begin{center}
\begin{tabular}{|c|c|c|}
\hline & \textit{natural increase} & \textit{natural mortalitity} \\
\hline
Malthus & $ \varepsilon N\left( t \right)$ & $ - \varepsilon N\left( t \right)$  \\
\hline \mbox{   } & \mbox{   } & \mbox{   }\\
Verhulst & $\varepsilon N\left( t \right)\left( {1 - \dfrac{N\left( t\right)}{K}} \right)$ & $ - \varepsilon N\left( t \right)\left( {1 + \dfrac{N\left( t \right)}{K}} \right)$ \\
\mbox{   } & \mbox{   } & \mbox{   }\\
\hline
\end{tabular}
\end{center}
\caption[Functional responses for natural increase and decrease.]{Functional responses for natural increase and decrease.}
\label{tab1}
\end{table}

\begin{figure}[htbp]
\centerline{\includegraphics[width = 8cm,height = 8cm]{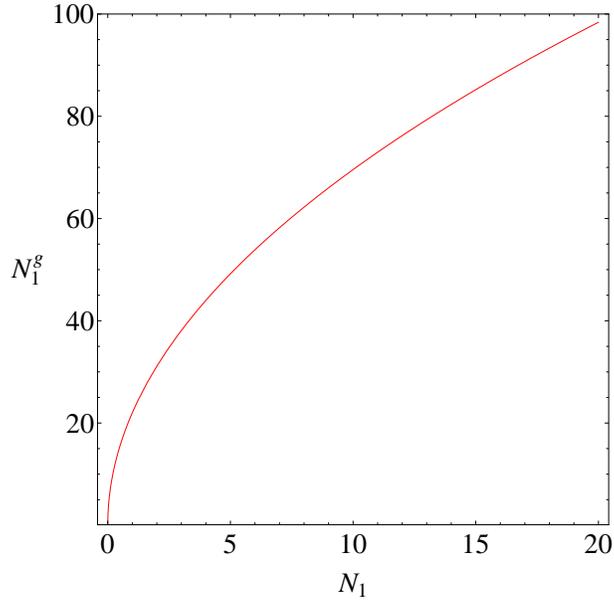}}
\caption{Gause's functional response}
\label{fig3}
\end{figure}

The functional response proposed by Volterra \cite{vol31} to describe predation and
which was based on the principle of encounters was proportional to the product
of the number of individuals of each species: $N_1\left( t \right) N_2 \left( t \right)$. In other words, the
predation rate was a "linear function" of the prey, that is, $N_1\left( t \right)$. A few years later,
Gause \cite{gau, gau2}, who was one of the first to make ``experimental verifications of the
mathematical theory of the struggle for life'', proposed another type of functional
response to describe predation, a ``nonlinear'' one, aimed at transcribing a certain
``satiety'' of the predator with respect to its prey: $N_1^g\left( t \right) N_2\left( t \right)$ with $0 < g \leqslant 1$. In
this case, the predation rate becomes indeed a ``nonlinear function'' of the prey, that
is,  $N_1^g\left( t \right)$ (see Fig. \ref{fig3}).

\begin{figure}[htbp]
\centerline{\includegraphics[width = 8cm,height = 8cm]{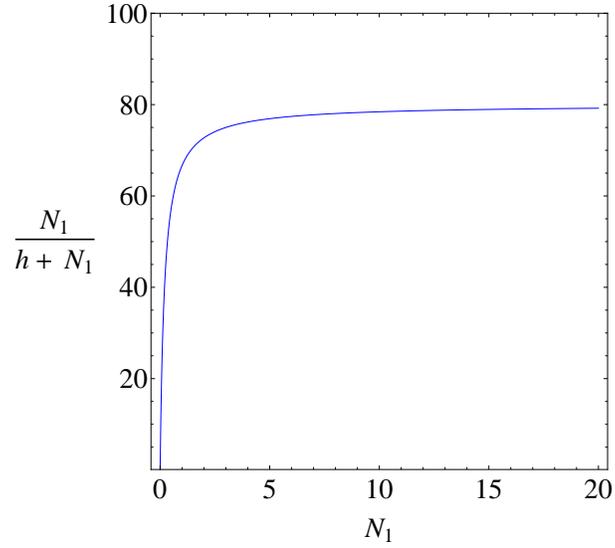}}
\caption{Functional response of Holling's type II.}
\label{fig4}
\end{figure}

In the late 1950s, entomologist Crawford Stanley Holling \cite{holla, hollb} developed, from
the celebrated ``disc equation'' two new functional responses for predation, also intended
to describe a certain satiety of the predator with respect to its prey: Holling
function of type II and Holling function of type III. This formulation assumes that
the predator divides its time between two kinds of activities: the search for its
prey, and its capture, which includes the time spent hunting, killing, devouring and
digesting it. Holling's type II (see \ref{fig4}) is a functional response in which the predator's ``attack'' rate increases when the prey number is low and then becomes constant when the predator reaches satiety. In other words, the predator causes
maximum mortality at low prey densities. Thus, type II functional responses are
typical of predators specializing in attacking one or few prey. In this case the
mortality of the prey decreases with their density. Holling's type II is represented
by:

\begin{equation}
\label{eq4}
\frac{N_1\left( t \right)}{h + N_1\left( t \right)}N_2\left( t \right)
\end{equation}

where $h$ represents \textit{half-saturation}, that is, the value of the prey density $N_1\left( t \right) = h$ for
which the predation level reaches a value equal to half its maximum.

\begin{figure}[htbp]
\centerline{\includegraphics[width = 8cm,height = 8cm]{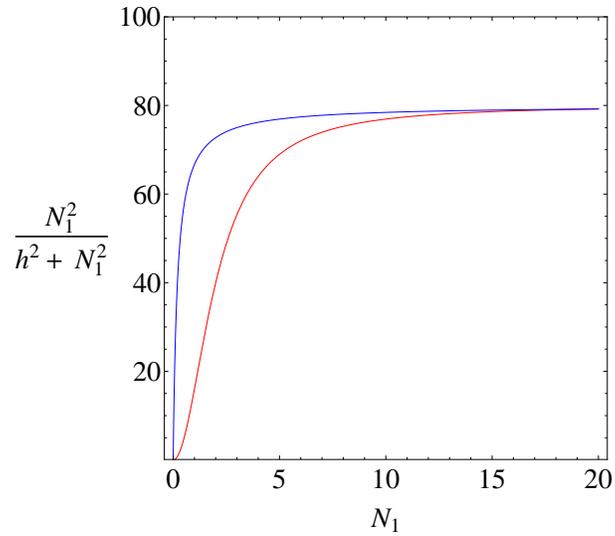}}
\caption{Functional response of Holling's type III (in red).}
\label{fig5}
\end{figure}

Holling's type III (see Fig. \ref{fig5} where Holling's type II is represented in blue and
Holling's type III in red) is a functional response in which the attack rate of the
predator first increases when the prey number is low and then slows down when the
predator reaches satiety. In other words, the predator increases its research activity
when prey density increases. Thus, type III functional responses are typical of
generalist predators moving from one species of prey to another and concentrating
their activities in areas where resources are abundant. In this case, the mortality
of the prey initially increases with their density and then decreases. Holling's type
III is represented by:

\begin{equation}
\label{eq5}
\frac{N_1^2\left( t \right) }{h^2 + N_1^2\left( t \right)}N_2 \left( t \right)
\end{equation}

where $h$ represents half-saturation, as above.

Contrary to what one might think, Holling did not deduce his two functional
responses from the observation of a natural environment, but, as he himself wrote
in his article, by concocting an artificial predatory-prey situation:

\begin{quote}

In the first artificial situation devised the ``prey'' were sandpaper
discs four centimetres in diameter thumb-tacked to a three-foot
square table. A blindfolded subject, the 'predator', stood in front
of the table and searched for the discs for one minute by tapping 
with her finger. As each disc was found, it was removed, set to
one side and searching continued. Each experiment was replicated
eight times at densities of discs ranging from four to 256 per nine
sq. ft.
The results of one such experiment are shown in Fig. 1 [see Fig. \ref{fig4}
here], where it can be seen that the number of discs picked up
increased at a progressively decreasing rate as the density of discs
rose. \cite[p. 385]{hollb}

\end{quote}

Mathematically, it has been proved by Real \cite{real} that the functional responses of
Holling type II and III are analogous to the function developed in 1913 by Leonor
Michaelis and Maud Menten \cite{mich} to describe the kinetics of enzymatic reactions.
All these functional responses are summarized Table 2.

\begin{table}[htbp]
\begin{center}
\begin{tabular}{|c|c|c|c|}

\hline \mbox{   } & \mbox{   } & \mbox{   } & \mbox{   }\\
Volterra & Gause & Holling type \textit{II} & Holling type \textit{III} \\
\mbox{   } & \mbox{   } & \mbox{   } & \mbox{   }\\
\hline \mbox{   } & \mbox{   } & \mbox{   } & \mbox{   }\\
$N_1 \left( t \right) N_2\left( t \right)$ & $N_1^g\left( t \right) N_2 \left( t \right)$ &
$\dfrac{N_1\left( t \right)}{h + N_1 \left( t \right)}N_2\left( t \right)$ &
$\dfrac{N_1^2\left( t \right)}{h^2 + N_1^2\left( t \right)}N_2\left( t \right)$ \\
\mbox{   } & \mbox{   } & \mbox{   } & \mbox{   }\\
\hline
\end{tabular}
\end{center}
\caption[Functional responses for predation]{Functional responses for predation}
\label{tab2}
\end{table}

\section{The origin of the paradox}

According to Yuri A. Kuznetsov \cite{kuz}, in the context of a Volterra predator-prey
model, the functional responses that limit growth (see Table 1) have a stabilizing
effect for the prey whereas the functional responses that limit predation (see Table
2) have a destabilizing effect for predators. Thus, in the Volterra predator-prey
model \ref{eq3}, by restricting the growth of the prey $N_1(t)$ by a Verhulst functional
response and the growth of the predator $N_2(t)$ by a functional response of Holling
type II, we obtain the celebrated Rosenzweig-MacArthur model \cite{rosen}:

\begin{equation}
\label{eq6}
\left\{
\begin{aligned}
\dfrac{dN_1}{dt} & = \varepsilon_1 N_1 - \lambda N_1^2 - \gamma_1 \frac{N_1}{h + N_1}N_2 \hfill \\
\dfrac{dN_2}{dt} & = - \varepsilon_2 N_2 + \gamma_2 \frac{N_1}{h + N_1}N_2
\end{aligned}
\right.
\end{equation}

where $\varepsilon_1$ represents the prey's growth rate in the absence of the predator,
$\gamma_1$ the predator's predation rate on the prey, $\varepsilon_2$ the predator's mortality rate in the absence
of prey, and $\gamma_2$ the predator's growth rate due to its predation. The combination of
these two effects then leads to the existence of a periodic solution. Henri Poincaré
\cite{poin1882} named this solution, or more exactly this periodic oscillation, \textit{limit cycle}. This
terminology derives from the fact that in the phase space ($N_1,N_2$), in this case the
prey-predator space, it takes the form of a cycle towards which every solution converges
asymptotically. Thus, the solution of the Rosenzweig-MacArthur model (\ref{eq6})
is periodic, like that of the Volterra model (\ref{eq3}). Nevertheless, there is a fundamental
difference in that it is totally independent of the initial conditions, that is, of the
initial density of prey and predators, which is consistent with reality. Indeed, in nature, the amplitude and period of the periodic oscillations of prey and predators cannot depend on their initial densities considered at a given time arbitrarily taken
as the origin of the time \cite{gin11,gin15,gin17}.

\section{Discussion}

Volterra's model, which in the study of nonlinear dynamic systems has become a
kind of paradigm, has been the subject of many books and articles criticizing it for
its lack of realism. The main flaws of this model are considered to be the absence
of limitation in the growth of prey and predator, the fact of not having taken into
account seasonality (that is, having considered the growth rate of prey and predator
as constant), and finally having proposed an ideal, simplified model by limiting it to
two species. The analysis of these flaws will shed light on the ``paradox'' of Volterra's
model. In the first place, it should be recalled that even today it is impossible to
formulate a simple expression of the periodic solutions of Volterra's model (\ref{eq3}), that
is, an expression that uses only elementary functions. Consequently, criticizing
Volterra for having proposed an ideal, simplified model is inadmissible. Indeed, it
is precisely thanks to the simplicity of his model that Volterra was able to solve
the problem posed by D'Ancona and explain the phenomenon by establishing the
law of perturbation of means. In fact, Volterra was well aware of the limitations
of the model and the assumptions attached to it, including the one assuming the
environment as invariable and without influence on the growth rate and the one that
considers the homogeneity of the individuals of each species. Indeed, in \cite{vol31} Volterra
laid the foundations for what would be much later called ``population dynamics''
and considered all the aspects of the problem: the influence of environment or
pollution, the heterogeneity of individuals or age groups (on hereditary actions \cite[p. 141]{vol31}), diffusion or migration (contribution of a small number of individuals \cite[p. 118]{vol31}), seasonality (variation of exterior conditions with time, \cite[p. 131]{vol31}). In
the introduction, he writes:

\begin{quote}
Certainly there exist periodic circumstances relating to environment,
as would be those, for example, which depend upon the changing
of the seasons, which produce forced oscillations of an external
character in the number of individuals of the various species.
These actions of external periodic nature were those which were
specially studied from the statistical point of view, but are there
others of internal character, having periods of their own which add
their action to these external causes and would exist even if these
were withdrawn? \cite[p. 5]{vol31}.
\end{quote}

Then, he adds:

\begin{quote}

Later, it will be observed that it is closer to reality to suppose that
the growth coefficients depend not only, at each moment, on the
current values of the quantities $N_i $ (characterizing the species $i$),
but also on past values up to a more or less remote period. It
will no longer be sufficient to consider them as functions of the
$N_i $, but as ``functionals'', and this will lead us to integro-differential
equations that we will approach from the equations we are led to in
the so-called ``hereditary'' mechanics \cite[p. 5]{vol31}.

\end{quote}

It is in this context that he invents the famous Volterra equations of the first
and second type. With regard to the limitation to a two-species model of predator-prey
type, in chapter II of \cite{vol31} Volterra presents the ``study of the coexistence of
any number of species''. Moreover, the chapter ends with the analysis of a ``very
remarkable case: that of three species of which the first feeds on the second and this
on the third''. Different types of interactions between species are also considered
outside of predation, including competition, cooperation and migration. Regarding
the absence of limitation in the growth of prey and predators, Volterra proposes
in chapter III of \cite{vol31} to study a model of $n$ coexisting species having reciprocal
actions, in which he replaces the Malthusian growth by a Verhulst logistic growth.
Thus, Italian or French authors who claim that Volterra's model is unrealistic have
obviously not read Volterra's texts \cite{vol27,vol28, vol31} and this is clearly one of the origins
of this paradox. Indeed, it seems that Volterra's model and criticisms of it are
known worldwide, but his work has not even been read by his detractors! For
authors who only read English, this paradox is rooted in another problem. While
the original Italian version of Volterra's study \cite{vol26} comprised 84 pages, the very
first English translation published in the prestigious journal Nature \cite{vol27} was only
two pages long\footnote{Under the same title, in a later issue of the same journal, two letters were published, one by
Lotka and one by Volterra himself \cite{lotkaN,vol27}.}. It is easy to understand that for many Volterra's study may have
seemed very limited given the conciseness of this summary, which did not allow
the exposition of analytical developments and considerably reduced the scope of
Volterra's results. Nevertheless, and this is again a paradox, even though by 1928
there was an integral English translation of the original Italian text by Volterra \cite{vol28},
it seems to have been almost totally ignored by English-speaking authors. Thus,
Volterra's work is perhaps the most quoted and the most criticized in the world
while being the least read and studied.

From 1926, Volterra's publications intensified in the field of population dynamics.
It was in 1931 that his work entitled \textit{Leçons sur la théorie mathématique de la lutte
pour la vie} \cite{vol31} was published in French\footnote{To the best of our knowledge, this book has never been translated into English.}, following a series of lectures at the Institut Henri Poincaré where Volterra had been invited by Borel. Compiled by Marcel
Brelot, this book contains Volterra's entire memoir on biological fluctuations \cite{vol26}
as well as a part about the case where heredity plays a role.
From 1936, Volterra questioned the validity of his model and the possibility of
experimental verification. The experiments of Gause (1910-1989) seemed to confirm
the first law, and he showed Volterra his acceptance of the model's forecasts. But
what followed was disappointing because they could not find irrefutable cases of a
cyclic behaviour in predator-prey ecosystems. Moreover, voiced by biologists Karl
Pearson (1857-1936) and Friedrich Simon Bodenheimer (1897-1959), objections
multiplied to the point of questioning the interpretation of D'Ancona's statistics,
which led to D'Ancona himself coming to doubt their merits. At present, the only
example of a predator-prey ecosystem showing a cyclic evolution is the famous set of
statistics by Hudson's Bay Company about hares and lynxes in Canada. According
to Giorgio Israel, it is thus important to note that unlike Van der Pol's model \cite{van}:

\begin{quote}
Volterra's model is not deduced from an analogy but from a more
traditional approach that consists in starting from the analysis of a real phenomenon, making some abstraction of the accessory aspects
such as friction, determining the state variables, formulate a
mathematical hypothesis concerning the pace of the phenomenon.
\ldots For Volterra's model, the problem of the experimental verification
consists in finding empirical evidence \textit{directly} verifying the laws
deduced from the model and not justifying its validity in an indirect
way, that is, from the effectiveness of some of its consequences.
Thus, for Volterra, the justification of D'Ancona's hypothesis on
the effects of fishing was not sufficient to prove the empirical validity
of the model. It is because of this conviction that he sought
for the rest of his life a direct \textit{empirical proof} of the validity of the
first law according to which the evolution of populations presents
periodic oscillations \cite{van}.

\end{quote}

Thus, if Volterra sought to make his mathematical research work at the service of
biological sciences, it seems that it always was with the aim of explaining phenomena
and describing reality as faithfully as possible:

\begin{quote}

`` \ldots the hypotheses will be seen to become more and more complex
in order to get closer to reality \ldots ''

\end{quote}

Translated from the French by Daniele A. Gewurz

\section*{Acknowledgments}

I wish to express my sincerest thanks to my friend Christian Gérini, Agrégé de
Mathématiques and Docteur en Histoire des Sciences, who encouraged and supported
me in my work.

\end{document}